\documentclass{article}
\usepackage{amsmath, amsfonts, amssymb}
\usepackage{natbib}
\usepackage{hyperref}
\newcommand{\Expect}[1]{\mathbb{E}\left[ #1 \right]}
\newcommand{\Var}[1]{\mathrm{Var}\left[ #1 \right]}
\newcommand{\Cov}[1]{\mathrm{Cov}\left[ #1 \right]}
\newcommand{\Prob}[1]{\mathbb{P}\left( #1 \right)}
\newcommand{\Indicator}[1]{\mathbb{I}\left( #1 \right)}
\newtheorem{lemma}{Lemma}
\newtheorem{theorem}{Theorem}
\newtheorem{corollary}{Corollary}
\newtheorem{proposition}{Proposition}
\newtheorem{definition}{Definition}
\begin{document}
\title{A Simple Non-Stationary Mean Ergodic Theorem, with Bonus Weak Law of
  Large Numbers}
\author{Cosma Rohilla Shalizi\thanks{Departments of Statistics and of Machine
    Learning, Carnegie Mellon University, 5000 Forbes Avenue, Pittsburgh, PA 15213, and the Santa Fe Institute, 1399 Hyde Park Road, Santa Fe, NM 87501.}}
\date{1 November 2021, small revision 19 March 2022}
\maketitle

\begin{abstract}
  This brief pedagogical note re-proves a simple theorem on the convergence, in
  $L_2$ and in probability, of time averages of non-stationary time series to
  the mean of expectation values.  The basic condition is that the sum of
  covariances grows sub-quadratically with the length of the time series.  I
  make no claim to originality; the result is widely, but unevenly, spread bit
  of folklore among users of applied probability.  The goal of this note
  is merely to even out that distribution.
\end{abstract}

I assume some familiarity with basic probability and stochastic processes,
along the lines of \citet{Grimmett-Stirzaker}, but re-prove a number of basic
results, to show that everything really is elementary.

Let $X_1, X_2, \ldots X_t, \ldots$ be a sequence of real-valued random
variables.  Assume that $\mu_t \equiv \Expect{X_t}$ and $\Cov{X_t,
  X_s}$ exist and are finite for all $t, s$.  The time
average of the $X_t$s,
\begin{equation}
A_n = \frac{1}{n}\sum_{t=1}^{n}{X_t} ~,
\end{equation}
converges on the average of the expectation values,
\begin{equation}
m_n = \frac{1}{n}\sum_{t=1}^{n}{\mu_t} ~,
\end{equation}
under a condition on the sum of the covariances,
\begin{equation}
V_n = \sum_{t=1}^{n}{\sum_{s=1}^{n}{\Cov{X_t, X_s}}} ~.
\end{equation}
The condition will be obvious after the following lemma, which will also be
useful for extensions.

\begin{lemma}
\begin{equation}
\Expect{(A_n - m_n)^2} =  \Var{A_n} = \frac{1}{n^2}V_n
\end{equation}
\label{lemma:mse-is-var-is-vn}
\end{lemma}

\textsc{Proof:}   Since, for any $Z$,
$\Expect{Z^2} = (\Expect{Z})^2 + \Var{Z}$, we have
\begin{equation}
\Expect{(A_n - m_n)^2} = (\Expect{A_n - m_n})^2 + \Var{A_n} \label{eqn:l2-distance}
\end{equation}
By linearity of expectation,
\begin{equation}
\Expect{A_n} = m_n
\end{equation}
On the other hand,
\begin{eqnarray}
\Var{A_n} & = & \frac{1}{n^2}{\Var{\sum_{t=1}^{n}{X_t}}}\\
& = & \frac{1}{n^2}{\sum_{t=1}^{n}{\sum_{s=1}^{n}{\Cov{X_t, X_s}}}}\\
& = & \frac{1}{n^2}{V_n}
\end{eqnarray}
through repeated application of the identity $\Var{B+C} = \Var{B}+\Var{C} +
2\Cov{B,C}$, and the definition of $V_n$.  Substituting into Eq.\ \ref{eqn:l2-distance},
\begin{equation}
\Expect{(A_n - m_n)^2} = 0^2 + \Var{A_n} = \frac{1}{n^2}V_n
\end{equation}
as was to be shown. $\Box$

\begin{theorem}[Mean ergodic theorem]
\label{thm:mean-ergodic}
  If $V_n = o(n^2)$, then $A_n - m_n \rightarrow 0$ in $L_2$ as $n\rightarrow \infty$.
\end{theorem}

\textsc{Proof:} Convergence in $L_2$ means that $\Expect{(A_n - m_n)^2}
\rightarrow 0$, and, by the lemma, $\Expect{(A_n - m_n)^2} = V_n/n^2$.  By
the assumption of the theorem, $V_n = o(n^2)$, hence $V_n/n^2 = o(1)$, so we
have
\begin{equation}
\Expect{(A_n - m_n)^2} = o(1) \rightarrow 0
\end{equation}
Thus $A_n - m_n \rightarrow 0$ in $L_2$; this is a mean (or mean-square) ergodic
theorem. $\Box$

Since $V_n$ is a sum of $n^2$ terms, if the sum is to be $o(n^2)$,
most of those terms must be shrinking rapidly to zero as $n$ grows.  That is,
$\Cov{X_t, X_{t+h}}$ must go to zero as $h\rightarrow\infty$.  Stationarity
(see immediately below) is not required.

\begin{definition}
The sequence of $X$s is {\bf weakly} or {\bf second-order stationary} when
$\Cov{X_t, X_{t+h}} = \gamma(h)$ for some
$\gamma$ and all $t$.
\end{definition}

If one {\em does} assume weak stationarity, a sufficient condition for $V_n =
o(n^2)$ is that $\sum_{h=-\infty}^{\infty}{\gamma(h)} = \tau \gamma(0) <
\infty$, since in that case $V_n = O(n)$.  $\tau$ has names like ``correlation
time'', ``autocorrelation time'', ``integrated autocovariance time'',
etc.\footnote{Some authors use these names for
  $\sum_{h=0}^{\infty}{\gamma(h)}$, or even for
  $\sum_{h=1}^{\infty}{\gamma(h)}$. Any one of these sums is finite if and only
  if the others are.}.  From Lemma \ref{lemma:mse-is-var-is-vn}, for
uncorrelated variables, $\Var{A_n} = \gamma(0)/n$, but when $0 < \tau <
\infty$, $n\Var{A_n} \rightarrow \tau\gamma(0)$, so the ``effective sample
size'' is reduced from $n$ to $n/\tau$.

It is unnecessary for the theorem to assume that $m_n$ converges, i.e., that
the instantaneous expectations $\mu_n$ are C{\`e}saro-convergent.  Still less
is it necessary to assume that the $\mu_n$ have a limit.  $\Var{X_t}$ need not
be constant or tending to a limit either, though it cannot grow too fast.

$L_2$ convergence implies convergence in probability, or what is usually
known as the weak law of large numbers.

\begin{corollary}[Weak law of large numbers]
Under the conditions of the theorem, $A_n - m_n \rightarrow 0$ in probability.
\end{corollary}

\textsc{Proof:} Use Chebyshev's inequality (re-proved below as Proposition
\ref{prop:chebyshev}): for any random variable $Z$,
\begin{equation}
\Prob{|Z - \Expect{Z}| \geq \epsilon} \leq \frac{\Var{Z}}{\epsilon^2}
\end{equation}
Applied to $A_n$, this gives, for each fixed $\epsilon$,
\begin{equation}
\Prob{|A_n - m_n| \geq \epsilon} \leq \frac{V_n n^{-2}}{\epsilon^2} \rightarrow 0
\end{equation}
which is the definition of convergence to 0 in probability.  $\Box$

The convergence results carry over easily to $d$-dimensional vectors, at least
if $d$ is fixed as $n$ grows.

\begin{corollary}
  Let $\vec{Y}_t = (Y_{t1}, Y_{t2}, \ldots Y_{td})$ be a sequence of
  $d$-dimensional vectors, with expected values $\vec{\mu}_t$, and define
  $V_{nj} = \sum_{t=1}^{n}{\sum_{s=1}^{n}{\Cov{Y_{tj}, Y_{sj}}}}$.  If $V_{nj}
  = o(n^2)$ for all $j$, then
  \begin{equation}
    \left\| \frac{1}{n}\sum_{t=1}^{n}{\vec{Y}_t} - \frac{1}{n}\sum_{t=1}^{n}{\vec{\mu}_t} \right\| \rightarrow 0
  \end{equation}
  in $L_2$ and in probability.
\end{corollary}
\textsc{Proof:} Apply the theorem and the previous corollary to each coordinate
of $Y$ separately to get convergence along each coordinate, and hence
convergence of the Euclidean distance to zero. $\Box$

If $V_n$ grows too fast, then the convergence {\em to a deterministic limit}
fails.

\begin{corollary}
If $V_n = \Omega(n^2)$, then $A_n  - m_n \not \rightarrow 0$ in $L_2$.
\end{corollary}
\textsc{Proof:} $V_n = \Omega(n^2)$ means that $\liminf{V_n/n^2} = v > 0$.
From the lemma, we know that $\Var{A_n} = V_n/n^2$, which does not go to zero,
so convergence in $L_2$ must fail.  $\Box$

Convergence {\em in probability} is slightly more delicate.

\begin{corollary}
If $V_n = \Omega(n^2)$ and $\Var{A_n^2} = O(V_n^2/n^4)$, then $A_n - m_n \not \rightarrow 0$ in probability.
\label{cor:no-conv-in-prob}
\end{corollary}
\textsc{Proof:} Begin with the previous corollary, and apply the 
Paley-Zygmund inequality (Proposition \ref{prop:paley-zygmund}) to the
non-negative random variable $(A_n - \mu_n)^2$, whose expected
value is (again, from Lemma \ref{lemma:mse-is-var-is-vn}) $V_n/n^2$.  By the
inequality, for any $\epsilon \leq V_n/n^2$,
\begin{eqnarray}
\Prob{(A_n - \mu_n)^2 \geq \epsilon} & \geq & \frac{(V_n/n^2 - \epsilon)^2}{\Var{(A_n - \mu_n)^2} + (V_n/n^2)^2}\\
\Prob{\left|A_n - \mu_n\right| \geq \sqrt{\epsilon}} & \geq & \frac{(V_n/n^2 - \epsilon)^2}{\Var{A_n^2} + (V_n/n^2)^2}
\end{eqnarray}
Restrict ourselves to $\epsilon < v/2$.  Then, for all sufficiently large $n$,
\begin{eqnarray}
\Prob{\left|A_n - \mu_n\right| \geq \sqrt{\epsilon}} & \geq & \frac{(v-\epsilon)^2}{v^2 + \Var{A_n^2}}\\
& \geq & \frac{(v/2)^2}{v^2 + \Var{A_n^2}} \\
& = & \frac{1}{4}\frac{1}{1+\Var{A_n^2}/v^2} > 0
\end{eqnarray}
so $A_n - \mu_n \not\rightarrow 0$ in probability.  $\Box$

{\em Remark 3:} I suspect the extra condition needed to force non-convergence
in probability can be weakened, because the underlying Paley-Zygmund inequality
used in the proof isn't necessarily sharp.  But an example helps show that {\em some} condition
is necessary.  Suppose that for each $t$, $X_t = \pm t^{3/2}$ with probability
$\frac{1}{2}t^{-2}$, otherwise $X_t = 0$, and that the $X_t$ are all mutually
independent.  Then $m_n=0$ for all $n$.  Moreover, $\Prob{X_t \neq 0} = 1/t^2$.
Since those probabilities are summable, by the Borel-Cantelli lemma, $\Prob{X_t
  \neq 0 ~\text{infinitely often}} = 0$.  But then $X_t =0$ for all but
finitely many $t$ almost surely, hence $A_n \rightarrow 0$ in probability.  On
the other hand, $\Var{X_t} = t$, so $V_n = n(n+1)/2$, $\lim{V_n/n^2} = 1/2$,
and $A_n \not\rightarrow 0$ in $L_2$.  Verifying that the second condition of
the corollary does {\em not} hold involves some straightforward but detailed
algebra, given in Appendix \ref{sec:gory-algebraic-details}.  While this
example is deliberately stylized, it does get at what's needed to have
convergence in probability without convergence in $L_2$: the probability that
$|A_n - m_n| \geq \epsilon$ has to be going to zero, no matter how small we set
$\epsilon$, but when there {\em are} fluctuations in $A_n$ away from $m_n$, they need to
be getting larger and larger.  Whether this is a realistic concern or a
paranoid fear will depend on the application.

{\em Remark 4:} Convergence of the $A_n$ not to the deterministic $m_n$ but to
a random limit, as in the full mean-square ergodic theorem for weakly
stationary sequences (as given by, e.g., \citet[\S 9.5, theorem
3]{Grimmett-Stirzaker}) would seem to require more advanced tools.

\paragraph{Credit}

I do not know the history of Theorem \ref{thm:mean-ergodic}, but I want to
emphasize again that it is not original to me.  I learned it, without
attribution to any particular source or even a name, when studying statistical
mechanics in the physics department of the University of Wisconsin-Madison in
the mid-1990s.  A version of the argument which assumes weak (second-order)
stationarity of the $X_t$ but allows for continuous time appears in \citet[pp.\
50--51]{Frisch-turbulence}.  The oldest version of the result I have been able
to locate is \citet{Taylor-diffusion-by-continuous-movements}.  (While the
paper was not published until 1922, it was read before the London Mathematical
Society in 1920.)  This again develops the result assuming weak stationarity,
but in both discrete and continuous time.  Taylor presents this as a new
result, but someone else might be able to claim historical priority.

\paragraph*{Acknowledgments}

I am grateful to David Darmon and Paul J. Wolfson for correspondence which led
me to write this; to support for a sabbatical year in 2017--2018 from Carnegie
Mellon University; and to my students in 36-462, ``Data over Space and Time'',
in 2018 and 2020, for letting me test versions of this material on them.

\bibliographystyle{crs}
\bibliography{locusts}

\begin{thebibliography}{3}
\newcommand{\enquote}[1]{``#1''}
\expandafter\ifx\csname natexlab\endcsname\relax\def\natexlab#1{#1}\fi
\expandafter\ifx\csname url\endcsname\relax
  \def\url#1{{\tt #1}}\fi
\expandafter\ifx\csname urlprefix\endcsname\relax\def\urlprefix{URL }\fi

\bibitem[\protect\citeauthoryear{Frisch}{Frisch}{1995}]{Frisch-turbulence}
Frisch, Uriel (1995).
\newblock {\em Turbulence: The Legacy of A. N. {Kolmogorov}\/}.
\newblock Cambridge, England: Cambridge University Press.

\bibitem[\protect\citeauthoryear{Grimmett and Stirzaker}{Grimmett and
  Stirzaker}{1992}]{Grimmett-Stirzaker}
Grimmett, G.~R. and D.~R. Stirzaker (1992).
\newblock {\em Probability and Random Processes\/}.
\newblock Oxford: Oxford University Press, 2nd edn.

\bibitem[\protect\citeauthoryear{Taylor}{Taylor}{1922}]{Taylor-diffusion-by-continuous-movements}
Taylor, G.~I. (1922).
\newblock \enquote{Diffusion by Continuous Movements.}
\newblock {\em Proceedings of the London Mathematical Society\/}, {\bf 20}:
  196--212.
\newblock
  \href{http://dx.doi.org/10.1112/plms/s2-20.1.196}{doi:10.1112/plms/s2-20.1.196}.

\end{thebibliography}

\clearpage

\appendix

\section{Upper Bounds: Markov and Chebyshev}

Going from $L_2$ convergence to convergence in probability uses an inequality
which has come to be associated with the name of Chebyshev:

\begin{proposition}[Chebyshev inequality]
\label{prop:chebyshev}
For any real-valued random variable $Z$,
\[
\Prob{|Z - \Expect{Z}| \geq \epsilon} \leq \frac{\Var{Z}}{\epsilon^2}
\]
\end{proposition}

This is itself an easy consequence of another inequality:

\begin{proposition}[Markov inequality]
\label{prop:markov}
For any non-negative, real-value random variable $Z$,
\[
\Prob{Z \geq \epsilon} \leq \frac{\Expect{Z}}{\epsilon}
\]
\end{proposition}

The intuition behind Markov's inequality is simple: the probability of $Z$
being large can't be too big, without also driving up the expected value of
$Z$.

\textsc{Proof} (of Proposition \ref{prop:markov}): For any event
$B$, $\Indicator{B} + \Indicator{B^c} = 1$.  So, clearly,
\[
Z = Z\Indicator{Z \geq \epsilon} + Z\Indicator{Z < \epsilon}
\]
and thus
\begin{eqnarray}
\label{eqn:expect-of-z} \Expect{Z} & = & \Expect{Z \Indicator{Z \geq \epsilon}} + \Expect{Z \Indicator{Z < \epsilon}}\\ 
& \geq & \Expect{Z \Indicator{Z \geq \epsilon}}\\
& \geq & \Expect{\epsilon \Indicator{Z \geq \epsilon}}\\
& = & \epsilon\Expect{\Indicator{Z \geq \epsilon}} = \epsilon \Prob{Z \geq \epsilon}
\end{eqnarray}
as was to be shown.  $\Box$

\textsc{Proof} (of Proposition \ref{prop:chebyshev}):
$|Z - \Expect{Z}| \geq \epsilon$ if and only if
$(Z - \Expect{Z})^2 \geq \epsilon^2$.  But $(Z - \Expect{Z})^2$ is a
non-negative, real-valued random variable, with expected value $\Var{Z}$, so
the proposition follows by applying Proposition \ref{prop:markov}. $\Box$

\section{Lower Bounds: Paley-Zygmund}

Markov's inequality says that the probability of large values can't be too
high, without increasing the expected value.  A counter-part inequality
essentially says that the probability of large values can't be too small,
either, without decreasing the expected value.  Start with
Equation \ref{eqn:expect-of-z}, again assuming $Z \geq 0$:
\begin{eqnarray}
\Expect{Z} & = & \Expect{Z \Indicator{Z \geq \epsilon}} + \Expect{Z \Indicator{Z < \epsilon}}\\
& \leq &  \Expect{Z \Indicator{Z \geq \epsilon}} + \epsilon\Expect{\Indicator{Z < \epsilon}}\\
\label{eqn:use-cauchy-schwarz} & \leq & \sqrt{\Expect{Z^2}\Expect{\Indicator{Z \geq \epsilon}}} +\epsilon\Expect{\Indicator{Z < \epsilon}}\\
& \leq & \sqrt{\Expect{Z^2}\Prob{Z \geq \epsilon}} + \epsilon\\
\frac{(\Expect{Z} - \epsilon)^2}{\Expect{Z^2}} & \leq & \Prob{Z \geq \epsilon}
\end{eqnarray}
where Eq.\ \ref{eqn:use-cauchy-schwarz} uses the Cauchy-Schwarz inequality.
We have thus proved

\begin{proposition}[Paley-Zygmund Inequality]
\label{prop:paley-zygmund}
For a random variable $Z \geq 0$, and $\epsilon \leq \Expect{Z}$,
\begin{equation}
\Prob{Z \geq \epsilon} \geq \frac{(\Expect{Z} - \epsilon)^2}{\Var{Z} + \Expect{Z}^2}
\end{equation}
\end{proposition}

(The proposition is usually stated in the form
\begin{equation}
\Prob{Z \geq \theta \Expect{Z}} \geq \frac{(1-\theta)^2\Expect{Z}^2}{\Var{Z} + \Expect{Z}^2}
\end{equation}
for $\theta \in (0,1)$, but this is clearly equivalent.)

\clearpage

\section{Algebraic Details for Remark 3}
\label{sec:gory-algebraic-details}

The example in Remark 3 converges in probability but not in $L_2$, so we need
to verify that one or the other conditions of Corollary
\ref{cor:no-conv-in-prob} fails.  In fact, the failing condition is the one
about fourth moments, that $\Var{A_n^2} = O(V_n^2/n^4)$.

\begin{eqnarray}
A_n^2 & = & n^{-2}\left[ \sum_{t=1}^{n}{X_t^2} + 2\sum_{t=1}^{n-1}\sum_{t+1}^{n}{X_t X_s} \right]\\
\Var{A_n^2} & = & n^{-4}\left[ \sum(\Var{X_t^2}) \right. \\
\nonumber & & +4\sum_{t=1}^{n-1}{\sum_{s=t+1}^{n}{\Var{X_t X_s}}}\\
\nonumber & &  + 2\sum_{t=1}^{n-1}{\sum_{s=t+1}^{n}{\Cov{X_t^2, X_s^2}}} \\
\nonumber & & + 4 \sum_{t=1}^{n}\sum_{s=1}^{n-1}{\sum_{r=s+1}^{n}{\Cov{X_t^2, X_s X_r}}} \\
\nonumber & & \left. + 8 \sum_{t=1}^{n-1}{\sum_{s=t+1}^{n}{\sum_{(r,q) \neq (t,s)}{\Cov{X_t X_s, X_r X_q}}}} \right]
\end{eqnarray}

Take the terms appearing in the expression for $\Var{A_n^2}$ one at a time:
\begin{eqnarray}
\Var{X_t^2} &= & \Expect{X_t^4} - (\Expect{X_t^2})^2\\
& = & t^6/t^2 - (\Var{X_t} + \Expect{X_t})^2\\
& = & t^4 - (t+0)^2 = t^4 - t^2 = t^2(t^2-1)\\
\Var{X_s X_t} & = & \Expect{X_t^2 X_s^2} - \left(\Expect{X_t X_s}\right)^2 \\
& = & \Expect{X_t^2} \Expect{X_s^2} - (\Expect{X_t}\Expect{X_s})^2 \\
& = & ts - 0 = ts
\end{eqnarray}
using the fact that the expectation of the product of two independent
variables is the product of their expectations.
On the other hand, for $s\neq t$,
\begin{eqnarray}
\Cov{X_t^2, X_s^2} & = & 0
\end{eqnarray}
by independence of the $X_t$s.  Similarly, unless $s=t$ or $r=t$, $\Cov{X_t^2, X_s X_r} = 0$.  Since $s\neq r$, at most one of $s$ and $r$ could equal $t$, so,
without loss of generality, say it's $s=t$.  Then we have
\begin{eqnarray}
\Cov{X_t^2, X_t X_r} & = & \Expect{X_t^3 X_r} - \Expect{X_t^2}\Expect{X_t X_r}\\
& = & \Expect{X_t^3}\Expect{X_r}  - \Expect{X_t^2}\Expect{X_t}\Expect{X_r}\\
& = & 0
\end{eqnarray}
So $\Cov{X_t^2, X_t X_r} = 0$
as well.  Likewise $\Cov{X_t X_s, X_r X_q} = 0$ if all the indices are
distinct, but even if indices overlap, everything boils out to zero.

Going back to the variance of $A_n^2$, then,
\begin{eqnarray}
\Var{A_n^2} & = & \frac{1}{n^4}\left[\sum_{t=1}^{n}{t^4 - t^2} + \sum_{t=1}^{n-1}{\sum_{s=t+1}^{n}{ts}}\right]\\
\end{eqnarray}
Now $\sum{t^4} \sim n^5$, so $\Var{A_n^2} \sim n$, but $V_n^2/n^4=
(n(n+1)/2)^2/n^4 \rightarrow 1/4$, so the second condition of the corollary is
indeed violated.

\end{document}